\documentclass{gOMS2e}

\usepackage{epstopdf}% To incorporate .eps illustrations using PDFLaTeX, etc.
\usepackage{subfigure}% Support for small, `sub' figures and tables

\usepackage{color}
\theoremstyle{plain}% Theorem-like structures

\theoremstyle{definition}

\theoremstyle{remark}

\newcommand{\R}{\mathbb{R}}

\renewcommand{\vec}[1]{\boldsymbol{#1}}

\begin{document}

%\jvol{00} \jnum{00} \jyear{2015} \jmonth{July}

% \articletype{GUIDE}

\title{A Non-Intrusive Parallel-in-Time Approach for Simultaneous Optimization with Unsteady PDEs}

\author{
\name{S. G\"unther\textsuperscript{a}$^{\ast}$\thanks{$^\ast$Corresponding author. Email: stefanie.guenther@scicomp.uni-kl.de}
  and N.R. Gauger\textsuperscript{a}
  and J.B. Schroder\textsuperscript{b}$^{\dagger}$\thanks{$^{\dagger}$This work was performed under the auspices of the U.S. Department of Energy by Lawrence Livermore National Laboratory under Contract DE-AC52--07NA27344, LLNL-JRNL-744565}
  }
\affil{\textsuperscript{a}Chair for Scientific Computing, TU Kaiserslautern, Kaiserslautern, Germany\\
       \textsuperscript{b}Center for Applied Scientific Computing, Lawrence Livermore National Laboratory, Livermore, CA
  }
% % \received{v5.0 released July 2015}
}

\maketitle

\begin{abstract}
  This paper presents a non-intrusive framework for integrating existing unsteady partial differential equation (PDE) solvers into a parallel-in-time simultaneous optimization algorithm. The time-parallelization is provided by the non-intrusive software library XBraid~\cite{xbraid-package},
  which applies an iterative multigrid reduction technique to the time domain of existing time-marching schemes for solving unsteady PDEs. Its general user-interface has been extended in \cite{guenther2017xbraid}
  for computing adjoint sensitivities such that gradients of output quantities with respect to design changes can be computed parallel-in-time alongside with the primal PDE solution. In this paper, the primal and adjoint XBraid iterations are embedded into a simultaneous optimization framework, namely the One-shot method. In this method, design updates towards optimality are employed after each state and adjoint update such that optimality and feasibility of the design and the PDE solution are reached simultaneously. The time-parallel optimization method is validated on an advection-dominated flow control problem which shows significant speedup over a classical time-serial optimization algorithm.

\end{abstract}


\begin{classcode}35Q93; 49M05; 93C20 \end{classcode}

\section{Introduction}
\label{sec:intro}

The aim of this paper is to develop a non-intrusive framework for reducing runtimes of conventional gradient-based optimization algorithms for optimal control and optimization with unsteady partial differential equations (PDEs).
Applications in aerodynamics include e.g. optimal active flow control, optimal shape design, inverse design or noise reduction \cite{nadarajah2007optimum,economon2013unsteady,rumpfkeil2010optimal,NeOeGaKrTh2016b,ZhAlGaICCFD2016,akcelik2002parallel}.

Despite the rapid increase of high-performance computing resources, computational runtimes associated with conventional optimization techniques are often tremendous and can easily scale up to weeks or even a month (compare, e.g., \cite{nemili2017accurate,nielsen2010discrete}).
% has stimulated the use of high-fidelity simulation tools for unsteady PDEs in an optimization and design framework.
The high computational costs associated with gradient-based optimization with unsteady PDEs can be traced back to the following factors.
First of all, the underlying unsteady dynamics need to be resolved properly. Existing simulation tools for unsteady PDEs typically involve a time integration loop that propagates the state of the described dynamical system forward in time, while solving nonlinear equations at each discrete time step. Hence, computational costs can be tremendous and grow linearly with the time-domain length.
 % even in the era of high-performance computing -- solving the underlying unsteady PDE still remains a computationally demanding task.
% The serial time propagation hence entails runtimes that increase linearly with the time-domain length.
Second of all, sensitivity information needs to be calculated in each iteration of a gradient-based optimization method. Here, the adjoint approach has become a powerful tool to compute the gradient and its computational costs are comparable to that of a pure simulation of the dynamics. When applied to unsteady PDEs, the adjoint approach involves a reverse time integration loop that propagates sensitivities backwards through the time domain \cite{nadarajah2007optimum,rumpfkeil2007general}. Hence, evaluating the gradient requires a forward loop in time to approximate the PDE solution followed by a backwards-in-time loop for the adjoint. If long time domains are considered, the serial forward and backward time propagation become a major bottleneck for fast and scalable optimization algorithms -- especially as speedups on future computer architectures will require more concurrency, as clock speeds remain stagnant.

Lastly, the choice of a suitable optimization technique often requires a careful trade-off between the computational efficiency against the intrusiveness with respect to the underlying PDE solver. Common approaches for unsteady PDE-constrained optimization are so-called \textit{reduced-space} optimization techniques that employ iterative design updates, while a corresponding PDE solution is recovered after each design change. Thus, the PDE-constraint is treated implicitly so that only minimal interactions of the optimizer and the PDE solver is required. However, the repeated forward and backward time integrations solving the PDE and the adjoint equations, respectively, entail an enormous runtime overhead for the optimization, when compared to a pure simulation. To address this issue, so-called \textit{full-space} or \textit{simultaneous} optimization approaches target reducing the optimization overhead by integrating the simulation directly into the optimization process, thus solving the unsteady PDE, the adjoint equations and the design optimization problem simultaneously. However, these approaches often require major modifications to the underlying simulation solver and are therefore far less commonly used for unsteady PDE-constrained optimization.

This paper tackles the above challenges in the following way: In order to reduce the runtime of the forward and the backward time loop, the unsteady PDE and corresponding adjoint equations are enhanced by a parallel-in-time integration approach that distributes the workload to multiple processors along the time domain. To this end, the non-intrusive software library XBraid \cite{xbraid-package} is employed, which adds a time-parallel capability to existing serial time-stepping codes. XBraid accesses the time-stepping routine of the original simulation code and applies an iterative multigrid reduction in time algorithm to the space-time problem. It then converges to the same solution as the existing serial simulation code, but can achieve a speedup due to new concurrency in the time domain.
In \cite{guenther2017xbraid}, XBraid's multigrid iterations have been enhanced with the ability to compute discrete adjoint sensitivities. Utilizing techniques from Automatic Differentiation (AD), the adjoint code runs backwards through the primal XBraid computations and accumulates sensitivity information parallel-in-time alongside the primal computations. Similar to the primal interface, existing adjoint time-marching schemes can be integrated through an extended adjoint user-interface such that a non-intrusive time-parallelization for existing adjoint codes is achieved.

While \cite{guenther2017xbraid} parallelized the state and adjoint solvers in time with XBraid, this paper novelly integrates that work into the context of an optimization process.  We choose the non-intrusive simultaneous \textit{One-shot} optimization method, because it reduces the optimization overhead of common reduced-space optimization techniques, allowing for greater speedups.
 The One-shot method has originally been developed for optimization with steady-state PDEs that are solved by fixed-point iterations \cite{bosse2014oneshot,schulz2004aerodynamic,gauger2009singlestep,oezkaya2014one}.
 These iterations are enhanced by adjoint and design updates such that design changes are realized already at an early stage of the simulation process.
In contrast to conventional reduced-space optimization techniques, the simultaneous One-shot approach  updates the design based on an approximate gradient obtained by solving for the PDE, the adjoint and the design optimization problem simultaneously.
% Since the primal fixed-point iterations of the PDE solver are treated in a black-box fashion, the simultaneous One-shot method provides a non-intrusive approach for reducing the runtime overhead of conventional reduced-space optimization techniques.

In this paper, the simultaneous One-shot method is applied to the time-parallel fixed-point iterations of XBraid's multigrid scheme. As those iterations update the PDE state over the entire time domain, they are well suited for integration into the non-intrusive simultaneous One-shot framework. In this setting, design updates are performed based on the approximate gradient that is available after each state and adjoint multigrid iteration.

Introducing time-parallelization to optimization and optimal control of unsteady PDEs is currently under active development. Recent approaches include, e.g., time-domain decomposition methods as in \cite{deng2016parallel,Heinkenschloss2005169,Ridzal2008}, Schwarz preconditioner approaches \cite{BaSt2015,Gander2016}, as well as preconditioners based on the parareal algorithm \cite{du2013inexact,ulbrich2015preconditioners}. In \cite{GoetschelMinion2017}, the parallel-in-time PFASST method is used to solve the primal and adjoint equation in an reduced-space optimization approach for parabolic optimal control problems. XBraid offers some advantages, such as multigrid scalability due to multiple time-levels as well as non-intrusive primal and adjoint user-interfaces, which allow for the re-use of existing simulation codes. Combined with the non-intrusive framework of full-space One-shot optimization, these strategies enable an easy transition from a conventional gradient-based optimization algorithm to parallel-in-time simultaneous optimization.

The paper is organized as follows: Section \ref{sec:unsteady_optim} introduces the optimization problem with unsteady PDE constraints. Conditions for optimality are displayed as well as a short discussion on conventional reduced-space optimization techniques which serve as a reference for comparing runtimes of the time-parallel One-shot method. The time-parallelization provided by XBraid is introduced for the forward- as well as backward-in-time state and adjoint simulations, respectively, in Section \ref{sec:xbraid}. In Section \ref{sec:oneshot}, the primal and adjoint XBraid iterations are integrated into the simultaneous One-shot optimization framework. Finally, numerical results are discussed in Section \ref{sec:numerics} which applies the scheme to an advection-dominated optimal control problem.

\section{Optimization with unsteady PDEs}
\label{sec:unsteady_optim}
We assume that a simulation code of time-marching type for solving the underlying unsteady PDE is available. Starting from a fixed initial condition $\vec u^0\in\R^M$, the time-marching simulation code propagates discretized states of the unsteady PDE forward with discrete time steps:
 \begin{align}
   \vec u^{i} &= \Phi^i(\vec u^{i-1}, \vec u^{i-2}, \dots, \rho) \quad \text{for} \quad i=1, \dots, N,
 \end{align}
where $\Phi^i\colon \R^M \times \R^p\to \R^M$ denotes a time integration scheme. It successively computes the discretized states $\vec u^i \in \R^M$ at discrete time step $0\leq t^i \leq T$ based on information at the previous time steps $\vec u^{i-1}, \vec u^{i-2}, \dots$, as well as certain design parameters $\rho \in\R^p$ that determine the PDE state. As implicit time integration methods are often preferred due to better stability properties, one application of the time-stepper $\Phi^i$ will likely require solving nonlinear and linear equations iteratively at each time step, as for example in a dual time-stepping framework \cite{jameson1991dualtime}. In order to facilitate readability, the dependency on previous time steps other than $\vec u^{i-1}$ will be dropped for the rest of this paper.

To transition from simulation to optimization, an objective function $J:\R^{N\times M}\times \R^p \to \R$ is defined that measures the time-average of some instantaneous quantity of the unsteady dynamics:
  \begin{align}
     J(\vec u, \rho) = \frac 1N \sum_{i=1}^N f(\vec u^i, \rho), %\approx \frac 1T\int_0^T f(u(t), \rho) \, \mathrm{d}t
  \end{align}
 where $\vec u := (\vec u^1, \dots, \vec u^N) \in \R^{N\times M}$ and $f\colon \R^M\times \R^p \to \R$.
The optimization task is then to find an optimal combination of design $\rho$ and state $\vec u$ that minimize $J$:
\begin{align}\label{optimproblem}
  \begin{split}
    \min_{\vec u, \rho} J(\vec u,& \rho) \\
    \text{s.t.} \quad \vec u^i = \Phi^i(&\vec u^{i-1}, \rho) \quad \forall \, i=1, \dots, N,
  \end{split}
\end{align}
using the initial condition $\vec u^0$.

% \subsection{Necessary optimality conditions}\label{KKTconditions}
The necessary optimality conditions for \eqref{optimproblem} can be derived using the corresponding Lagrangian function
\begin{align}
  L(\vec u, \vec {\bar u}, \rho) := J(\vec u, \rho) + \sum_{i=1}^N \left(\vec {\bar u}^i\right)^T\left( \Phi^i(\vec u^{i-1}, \rho) - \vec u^i \right),
\end{align}
where $\vec{\bar u} = \left(\vec {\bar u}^1,\ldots, \vec {\bar u}^N\right) \in \R^{N\times M}$ denote the so-called adjoint variables. At an optimal point of the optimization problem \eqref{optimproblem}, the derivative of the Lagrangian function with respect to the state, the adjoint and the design variables is zero which leads to the following set of equations (the KKT conditions \cite{nocedal2006numerical}):
\begin{align}
  \vec u^i &= \Phi^i(\vec u^{i-1}, \rho) & \forall \, i=1, \dots, N \label{KKTstate}\\
  \vec{\bar u^i} &= \nabla_{\vec u^i}J(\vec u, \rho) + \left(\partial_{\vec u^{i}}\Phi^{i+1}(\vec u^{i}, \rho)\right)^T\vec{\bar u^{i+1}} & \forall \, i=N, \dots, 1 \label{KKTadjoint}\\
    0 &= \nabla_{\rho}J(\vec u, \rho) + \sum_{i=1}^N  \left( \partial_{\rho} \Phi^{i}(\vec { u}^{i-1}, \rho) \right)^T{\vec {\bar u}}^i, \label{KKTdesign}
\end{align}
where subscripts denote partial derivatives. Here, $\vec u^0$ is the fixed initial condition of the state equations \eqref{KKTstate} and $\vec {\bar u}^{N+1}:=0$ serves as a terminal condition for the adjoint equations in \eqref{KKTadjoint}. Starting from this terminal condition, the adjoint equations can be solved successively moving backwards in time from $t^N,\dots,t^1$, given a fixed design parameter $\rho$ and the corresponding PDE state $\vec u$. The adjoint variables $\vec{\bar u}$ can then be used to calculate the right hand side of the so-called design equation in \eqref{KKTdesign}. This right hand side reflects the total derivative of the objective function $J$ with respect to design changes and is therefore often referred to as the \textit{reduced gradient}. Hence, the adjoint variables can be interpreted as an auxiliary construct to calculate the sensitivity of the objective function with respect to design changes. This sensitivity is used to modify the design parameter $\rho$ in an outer optimization cycle.

% \subsection{sec:reducedspaceoptim}
Conventional reduced-space optimization algorithms typically utilize a preconditioned version of the reduced gradient in order to update the design variable. To that end, a preconditioning matrix $B_k^{-1}$ is chosen in each optimization iteration that approximates the Hessian of the objective function. A common approach is to approximate it using low rank update schemes which are based on the reduced gradient, e.g. BFGS updates \cite{nocedal2006numerical}. A conventional reduced-space optimization algorithm then reads as follows:
\begin{align}
  % \text{Initialize } \rho_0;
  \text{Iterate } k = &0,1,\dots \notag\\
  \text{For } i&=1,\dots,N: \,
       \vec u^i = \Phi^i(\vec u^{i-1}, \rho_k) \label{reducedoptimiter1}\\
  \text{For } i&=N,\dots,1: \,
      \vec{\bar u^i} = \nabla_{\vec u^i}J(\vec u, \rho_k) + \left(\partial_{\vec u^{i}}\Phi^{i+1}(\vec u^{i}, \rho_k)\right)^T\vec{\bar u^{i+1}} \label{reducedoptimiter2}\\
  \rho_{k+1} &= \rho_k - B^{-1}_k \left(\nabla_{\rho}J(\vec u, \rho_k) + \sum_{i=1}^N  \left( \partial_{\rho} \Phi^{i}(\vec { u}^{i-1}, \rho_k) \right)^T{\vec {\bar u}}^i\right),\label{reducedoptimiter3}
\end{align}
starting from an initial design $\rho_0$ and utilizing the initial condition $\vec u^0$ and adjoint terminal condition $\vec {\bar u}^{N+1}=0$. Each outer optimization cycle with the current design $\rho_k$ first computes the PDE state in a forward-in-time loop, followed by a backwards-in-time loop to calculate the adjoint variables. Then, the preconditioned reduced gradient is used to update the design variable in \eqref{reducedoptimiter3}.
 % The optimization algorithms terminates, if the norm of the reduced gradient is below a prescribed tolerance such that \eqref{KKTdesign} holds up to a certain accuracy.
 This algorithm serves as a reference for computing runtime speedups of the proposed time-parallel One-shot optimization method.

%
% \begin{algorithm}[h]
%   \caption{Conventional reduced-space optimization procedure}\label{alg:reducedoptim}
%   \begin{algorithmic}[0]
%     \State Initialize $\rho_0$
%     \State \textbf{for} {$k=0,1,2,\dots$} \textbf{do}
%     \State $\quad$ 1. Solve the unsteady PDE:
%     \State $\qquad \,$ \textbf{for} {$i=1,\dots,N$} \textbf{do} %\Comment{Solve the unsteady PDE}
%          \Statex $\qquad \qquad  \vec u^i = \Phi^i(\vec u^{i-1}, \rho_k)$
%     \State $\quad$ 2. Solve the unsteady adjoint equations:
%     \State $\qquad \,$ \textbf{for }{$i=N,\dots,1$} \textbf{do} %\Comment{Solve the unsteady adjoint equations}
%       \State $\qquad \qquad \vec{\bar u^i} = \nabla_{\vec u^i}J(\vec u, \rho_k) + \left(\partial_{\vec u^{i}}\Phi^{i+1}(\vec u^{i}, \rho_k)\right)^T\vec{\bar u^{i+1}}$
%
%     \State $\quad$ 3. Update the design:
%     % \State $\qquad \,$ $\bar\rho = \nabla_{\rho}J(\vec u, \rho_k) + \sum_{i=1}^N  \left( \partial_{\rho} \Phi^{i}(\vec { u}^{i-1}, \rho_k) \right)^T{\vec {\bar u}}^i$
%     \State $\qquad \,$ $\rho_{k+1} = \rho_k - B^{-1}_k \left(\nabla_{\rho}J(\vec u, \rho_k) + \sum_{i=1}^N  \left( \partial_{\rho} \Phi^{i}(\vec { u}^{i-1}, \rho_k) \right)^T{\vec {\bar u}}^i\right)$
%     \State \textbf{until} $\|\bar \rho \| \leq \epsilon$
%   \end{algorithmic}
% \end{algorithm}

\section{Non-intrusive time-parallelization}
\label{sec:xbraid}
This section introduces the time-parallelization for the forward- and backward time-marching loop of the PDE state and the adjoint variables, respectively. To this end, the parallel-in-time software library XBraid \cite{xbraid-package} is utilized, which applies a multigrid reduction in time (MGRIT) algorithm to the forward time-stepping scheme. Time-parallelization of the backwards-in-time adjoint loop is then described utilizing techniques from AD applied to MGRIT.

\subsection{Multigrid reduction in time}

% The multigrid reduction in time (MGRIT) algorithm is used to solve
% equations (\ref{reducedoptimiter1}) and (\ref{reducedoptimiter2}).
We consider first the forward problem (\ref{reducedoptimiter1}), in the linear case
for simplicity.  This problem is equivalent to solving
the block lower bidiagonal system defined by
\begin{equation}\label{eqn:linear-system}
   A \vec u :=
   \left(
   \begin{array}{cccc}
   I                &        &           & \\
   -\Phi^{1}        & I      &           & \\
                    & \ddots & \ddots    & \\
                    &        & -\Phi^{N} & I
   \end{array}
   \right)
   \left(
   \begin{array}{c}
    \vec u^{0} \\
    \vec u^{1} \\
    \vdots \\
    \vec u^{N}
   \end{array}
   \right), \mbox{ to yield }
   A \vec u =
   \left(
   \begin{array}{c}
    \vec u^{0} \\
    0 \\
    \vdots \\
    0
   \end{array}
   \right).
\end{equation}
Sequential time stepping is simply an $O(N)$ forward sequential solve of
(\ref{eqn:linear-system}).  MGRIT instead solves (\ref{eqn:linear-system}) with
an iterative, $O(N)$, and highly parallel multigrid reduction method
\cite{oosterlee_book, ries_etal_1983_mgr}.

In practice, the application of $\Phi^{i}$ also encapsulates any forcing terms,
making it an affine operator.  Thus in this setting, and in the more general
nonlinear setting, we use the Full Approximation Storage (FAS) nonlinear
multigrid cycling \cite{Brandt_1977} to solve the nonlinear system defined by
\begin{align}
  A (\vec u, \rho) := \begin{pmatrix}
      \vec u^0 \\
      \Phi^1(\vec u^0, \rho) - \vec u^1 \\
      \Phi^2(\vec u^1, \rho) - \vec u^2 \\
      \vdots \\
      \Phi^N(\vec u^{N-1}, \rho) - \vec u^N
    \end{pmatrix},
     \mbox{ to yield }
  A (\vec u, \rho) =
   \left(
   \begin{array}{c}
    \vec u^{0} \\
    0 \\
    \vdots \\
    0
   \end{array}
   \right).
\end{align}

\begin{figure}
\center
\includegraphics{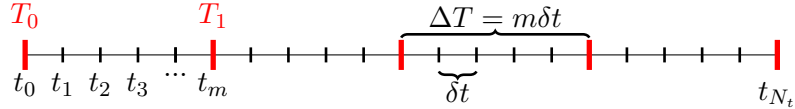}
% \usetikzlibrary{decorations.pathreplacing}
% \begin{tikzpicture}
% \path[draw] (0,0) -- (10,0);
%
% % Draw Fine Points
% \foreach \x in {0,1,...,3}
%    \path[draw, line width = 1pt] (\x/2,0.1) -- (\x/2,-0.1) node [below] {$t_{\x}$};
% \path[draw, line width = 1pt] (2,0.1) -- (2,-0.1) node [below] {$...$};
% \path[draw, line width = 1pt] (2.5,0.1) -- (2.5,-0.1) node [below] {$t_m$};
% \foreach \x in {4,5,...,20}
%    \path[draw, line width = 1pt] (\x/2,0.1) -- (\x/2,-0.1);
%
% % Draw Coarse Points
% \path[draw, line width = 2pt, color = red] (0,-0.2) -- (0,0.2) node [above] {$T_0$};
% \path[draw, line width = 2pt, color = red] (2.5,-0.2) -- (2.5,0.2) node [above] {$T_1$};
% \path[draw, line width = 2pt, color = red] (5,-0.2) -- (5,0.2) ;
% \path[draw, line width = 2pt, color = red] (7.5,-0.2) -- (7.5,0.2) ;
% \path[draw, line width = 2pt, color = red] (10,0.2) -- (10,-0.2) node [below] {$\textcolor{black}{t_{N_t}}$};
% \draw [decoration={brace}, decorate, line width = 1pt] (5,0.23) -- (7.5,0.23) node [above, pos=0.5] {$\Delta T = m \delta t$};
% \draw [decoration={brace,mirror}, decorate, line width = 1pt] (5.5,-0.22) -- (6,-0.22) node [below, pos=0.5] {$\delta t$};
% \end{tikzpicture}
\caption{Fine grid ($t_i$) and coarse grid ($T_j$) for coarsening factor
   $m=5$.  The fine points (black) are eliminated to yield a coarse level
   (red).}
\label{fig-grid}
\end{figure}

For a full description of the FAS-based MGRIT algorithm, see
\cite{schroder2015multigrid, falgout2014compressible}.  As a summary, the MGRIT method
is derived from approximate block cyclic reduction methods, which are
well-known methods for tridiagonal systems.  The time points are decomposed
into red and black points, as depicted in Figure \ref{fig-grid}.  Then, a
coarser time level is constructed by approximately eliminating the block rows
and columns in (\ref{eqn:linear-system}) that correspond to the black points.
The result is a smaller problem with the same form as
(\ref{eqn:linear-system}), but defined only at the red points.  The cyclic
reduction is approximate because the ideal coarse level time-stepping operator
is approximated by using the fine level time stepping scheme $\Phi^{i}$, but rediscretized
with the larger coarse level time step size.  This is done recursively, to
construct a multilevel hierarchy, where the coarser temporal levels compute
error corrections that accelerate the solution on the fine grid equations in
(\ref{eqn:linear-system}).

Complementing the coarse grid error correction is the block Jacobi relaxation scheme
that alternates between the red and black time points.
F-relaxation refers to block Jacobi applied to the block rows corresponding to the
black time points, and C-relaxation refers similarly to block Jacobi over
the red time points.
FCF-relaxation is then three sweeps over the F-, C- and then F-points.
Movement between the levels is done by injection, and the order of the levels
visited is done using either standard multigrid V-cycles, or the more powerful
F-cycles \cite{oosterlee_book}.
%\textbf{I can insert the FAS algorithm, but I'm not sure that we need it.
%Just let me know}

% In an analogous fashion, an MGRIT adjoint solver for equation
% (\ref{reducedoptimiter2}) can also be constructed \cite{guenther2017xbraid},
% except that now the system moves backwards in time.  Here, the system
% corresponding to (\ref{eqn:linear-system}) is now upper block bidiagonal.  The
% construction of this adjoint solver was the main topic in
% \cite{guenther2017xbraid}, and we refer the reader there.  In that work,

We use the XBraid framework \cite{xbraid-package}, which implements MGRIT. XBraid is non-intrusive due to general user interface functions that allow for the wrapping of existing serial time-stepping codes that solve (\ref{reducedoptimiter1}).

In summary, MGRIT provides an equivalent solution when compared to sequential
time stepping on the finest grid, to within a halting tolerance.  Also, both
MGRIT and sequential time stepping are optimal and $O(N)$; however, the
constant is higher for MGRIT, which creates a crossover point.  This implies
that a certain number of processors are required to create a speedup with MGRIT
\cite{falgout2014parallel}, and this will be evident in our results.

%Can you give an even briefer summary?  Remove the picture?  Just say decompose
%into red and black points.  See the other paper?  Or just condense the
%language.

\subsection{MGRIT as a time-parallel fixed-point solver}

MGRIT is an iterative solution scheme, and
as such, when it converges, it is a fixed-point scheme.  Let
$\vec u_k = \left(\vec u^1_k, \dots, \vec u^N_k\right) \in \R^{N\times M}$ be
the space-time solution after $k$ MGRIT iterations, and
$H\colon \R^{N\times M}\times \R^p \to \R^{N\times M }$
represent the action of one MGRIT iteration, i.e., one multigrid in time cycle.  Then,
the forward in time solution process is defined by
\begin{align}
   \text{for } k=0,1,\ldots: \quad \vec u_{k+1} = H(\vec u_k, \rho).
\end{align}
In the limit as $k \rightarrow \infty$, if MGRIT is convergent (i.e., contractive),
\begin{align}
   \vec u = H(\vec u, \rho), \mbox{ where } \vec u^i = \Phi^i(\vec u^{i-1}, \rho) \quad \forall \, i=1,\dots, N.
\end{align}

MGRIT has been shown to be contractive,
\begin{align}
   \|\partial_{\vec{u}}H(\vec u, \rho) \|_2 \leq \eta < 1,
\end{align}
i.e., a fixed-point method in a variety of settings.  The work
\cite{friedhoff2015} shows contractivity for linear parabolic problems in a
two- and three-level setting.  In \cite{dobrev2016twolevel}, further exploration
is done showing contractivity for basic linear parabolic and hyperbolic
problems, although we note that hyperbolic problems require a more careful
parameter tuning and generally converge more slowly than for the parabolic case. The
recent paper \cite{HeNoRoSc2017} shows contractivity for two-level MGRIT and
 a linearized elasticity formulation.
% I can give more math here on what the two-grid error propagator looks like,
% and how it can be bounded, but that's probably not needed.
While these results are limited to the linear setting with two- or
three-levels, they have been shown experimentally to be indicative of
multilevel results, including for nonlinear problems,
\cite{dobrev2016twolevel, falgout2014parallel, HeNoRoSc2017,
schroder2015multigrid}.
 % Lastly, we can also say that these contractivity
% results extend to the backwards in time MGRIT adjoint solver, because it is essentially
% equivalent to the forward in time MGRIT solver.

\subsection{Time-parallel adjoint sensitivities}
In \cite{guenther2017xbraid}, the XBraid library, which implements MGRIT, has been extended for computing time-parallel consistent and discrete adjoint sensitivities of the objective function with respect to design changes.
Since changes in $\rho$ induce changes of the objective function $J$ both explicitly as well as implicitly through the state variable solving $\vec u = H(\vec u, \rho)$, the total derivative of $J$ can be expressed by
\begin{align}\label{reducedgradient}
    \nabla_{\rho} J(\vec u, \rho) + \left({\partial_{ \rho} H(\vec u, \rho)} \right)^T\vec{\bar u},
\end{align}
where $\bar{\vec u}$ solves the adjoint equation
\begin{align}\label{adjointequation}
   \vec{\bar u} = \nabla_{\vec u} J(\vec u, \rho) + \left({\partial_{\vec u} H(\vec u, \rho)}\right)^T\vec {\bar u} .
\end{align}

Exploiting techniques from AD, the primal MGRIT iterations have been enhanced with an adjoint iteration that runs backwards through the code collecting the desired partial derivatives in \eqref{reducedgradient} and \eqref{adjointequation}.
The resulting scheme computes time-parallel adjoint sensitivities alongside the primal time-parallel calculation of the unsteady dynamcis in the following so-called \textit{piggyback} iteration:
\begin{align}
  \text{For } k=0,1,\dots&: \notag \\
  \vec u_{k+1} &=  H(\vec u_k, \rho) \label{piggyback_itertion_primal} \\
   \bar{\vec u}_{k+1}  &= \nabla_{\vec u}J(\vec u_k, \rho) + \left( {\partial_{\vec u} H(\vec u_k, \rho)}\right)^T \bar{\vec u}_k,\label{piggyback_itertion_adjoint}
\end{align}
for a given design $\rho$.
Since $H$ is contractive, the piggyback iteration converges simultaneously to the solution of the unsteady PDE as well as the adjoint equations with convergence rate $\eta$. However, since the adjoint equation \eqref{adjointequation} depends on $\vec u$, it is expected that the adjoint iterations exhibits a certain time-lag when compared to the primal one which has been analyzed in \cite{griewank2002reduced}.

 After the piggyback iteration has converged, the adjoint variable $\vec{\bar u}$ satisfies
\begin{align}
   \bar{\vec u} &= \nabla_{\vec u}J(\vec u, \rho) + \left( {\partial_{\vec u} H(\vec u, \rho)}\right)^T \bar{\vec u} \\
   \Leftrightarrow \quad \vec {\bar u}^i &= \nabla_{\vec u^i} J(\vec u,\rho) + \left( \partial_{\vec u^i} \Phi^{i+1}(\vec u^i, \rho)\right)^T\vec{\bar u^{i+1}}
   \quad  \forall \, i=N,\dots,1.
 \end{align}
Thus, the adjoint equations at each time step defined in \eqref{KKTadjoint}, Section \ref{sec:unsteady_optim}, are satisfied.

Similar to the primal XBraid implementation, the XBraid adjoint solver is non-intrusive in the sense that existing adjoint simulation codes can easily be integrated through extended user interface routines that define an adjoint time step for solving
(\ref{reducedoptimiter2}).

\section{Simultaneous One-shot optimization}
\label{sec:oneshot}
In this section, the piggyback iteration \eqref{piggyback_itertion_primal}--\eqref{piggyback_itertion_adjoint} is integrated into a simultaneous One-shot optimization framework. To this end, a reformulated optimization problem is considered that utilizes MGRIT's fixed-point equation as constraints:

  \begin{align}\label{unsteadyoptimXBraid}
      \min_{\vec u, \rho} J(\vec u, \rho) \quad \text{s.t.} \quad \vec u = H(\vec u, \rho).
  \end{align}

  At any optimal point of \eqref{unsteadyoptimXBraid}, the state satisfies the constraint and the gradient of the objective function with respect to the design is zero. Hence, the following set of equations hold:
  \begin{align}
    \vec u &= H(\vec u, \rho) \\
    \bar{\vec u} &= \nabla_{\vec u} J(\vec u, \rho) + \left(\partial_{\vec u} H(\vec u, \rho) \right)^T \bar{\vec u} \\
    0 &= \nabla_{\rho} J(\vec u, \rho) + \left(\partial_{\rho} H(\vec u, \rho) \right)^T \bar{\vec u}.
  \end{align}
Instead of solving these systems exactly and in consecutive order, the simultaneous One-shot optimization method aims at iteratively solving the whole system in a coupled iteration for the state, the adjoint and the design variables. In particular, the following iteration is proposed to find an optimal point satisfying the above equations:
\begin{align}
  \text{Iterate } \, k=0, &1, \dots \notag\\
  \vec u_{k+1} &= H(\vec u_k, \rho_k) \label{oneshotiter_state}\\
  \bar{\vec u}_{k+1} &= \nabla_{\vec u} J(\vec u_k, \rho_k) + \left(\partial_{\vec u} H(\vec u_k, \rho_k) \right)^T \bar{\vec u}_k \label{oneshotiter_adjoint}\\
  \rho_{k+1} &= \rho_k - B^{-1}_k \left( \nabla_{\rho} J(\vec u_{k}, \rho_k) + \left(\partial_{\rho} H(\vec u_k, \rho_k) \right)^T \bar{\vec u}_k  \right), \label{oneshotiter_design}
\end{align}
starting from some initial guess on the design $\rho_0$, the state $\vec u_0$ and the adjoint variable $\vec{\bar u}_0$. The above One-shot iteration performs updates of all three variables in a coupled iteration based on their current approximations. Thus, instead of exactly solving for the state and the adjoint variables before each design change, the One-shot method performs design updates that are based on the inexact reduced gradient. The first two lines \eqref{oneshotiter_state} and \eqref{oneshotiter_adjoint} form one piggyback iteration which is executed in parallel across the time domains. The last line \eqref{oneshotiter_design}, however, involves communication over all processors in order to collect and add up the desired sensitivity and broadcast the new design to the time domains.

The preconditioning matrix $B_k^{-1}$ in \eqref{oneshotiter_design} needs to compensate for the inexactness of the unsteady dynamics during optimization in order to ensure convergence of the One-shot method. It is therefore proposed in \cite{hamdi2011reduced}, and with further simplification in \cite{BlDeBaGaRe2017}, to search for the descent of an augmented Lagrangian function
\begin{align}
  L^{a}(\vec u, \bar{\vec u}, \rho) &= \frac{\alpha}{2} \|H(\vec u, \rho) - \vec u \|^2 +  J(\vec u, \rho) + (\bar{\vec u})^T\left(H(\vec u, \rho) - \vec u\right),
\end{align}
where the weighted residual of the state equation has been added to the Lagrangian function with $\alpha>0$. If $\alpha$ is big enough, in particular if
\begin{align}
  \alpha > \frac{2l}{1-\eta},
\end{align}
% JBS Comment:  Is r defined anywhere?
% SG: There shouldn't be an r at all! I corrected the formula. Good that you saw that!
where $\eta<1$ is the contractivity rate of the MGRIT fixed-point iteration and $l$ quantifies the adjoint time lag with $\|\vec{\bar u}_{k+1} - \vec {\bar u}_k\| \leq l\|\vec u_{k+1} - \vec u_k \|$, then $L^a$ is an exact penalty function such that optimal points of the optimization problem coincide with minima of $L^a$. Further, it can be shown that if the preconditioner $B_k$ approximates the Hessian of $L^a$, each One-shot iteration \eqref{oneshotiter_state} -- \eqref{oneshotiter_design} yields descent on $L^a$. Hence, the One-shot method converges to an optimal point.

Consequently, a suitable preconditioner can be approximated numerically using low-rank update schemes (e.g. BFGS) that are based on the gradient of the augmented Lagrangian. Quantifying $\alpha$, however, might be somewhat critical as the contractivity rate $\eta$ may not be known in advance. It is therefore a common approach to use BFGS updates based on the reduced gradient, instead of the gradient of $L^a$. While this choice corresponds to setting $\alpha = 0$, convergence has already been observed numerically \cite{bosse2014oneshot}.

\section{Numerical Results}
\label{sec:numerics}

The time-parallel One-shot method is applied to an inverse design problem subject to an advection-dominated model problem. The model problem is chosen such that it mimics flow dynamics past cylindrical bluff bodies at low Reynolds numbers. The near wake behind bluff bodies is dominated by a recirculation zone where regular periodic vortices are forming. Those vortices then shed into the far wake where they slowly dissipate \cite{cylinder}. The two dynamical zones are mimicked by the following approaches. In the near wake, a nonlinear ODE exhibiting self-excited oscillations is utilized, namely the Van-der-Pol oscillator \cite{leweke1994model,kanamaru2007van}. The far wake is modeled by an advection-diffusion equation whose upstream boundary condition is determined by the oscillating ODE mimicking the near wake.

The modeling equations under consideration then read
\begin{align}
  \partial_t v(t,x) + a\partial_x v(t,x) - \mu \partial_{xx} v(t,x) &= 0 \qquad \forall x \in (0,1), t\in (0,T)   \\
  v(t,0) - \mu \partial_x v(t,0) &= z(t) \quad  \forall t\in(0,T) \\
  \partial_{xx}v(t,1) &= 0  \qquad  \forall t\in (0,T)  \\
  v(0,x) &= 1 \qquad \forall x \in [0,1],
\end{align}
where the advection term dominates with $a=1$ and a small diffusion term with $\mu=10^{-5}$ is added. The upstream boundary $z(t)$ is determined through the Van-der-Pol oscillator:
\begin{align}\label{PinT:Numerics_VDP}
  \begin{pmatrix}
    \partial_t z(t) \\ \partial_t w(t)
  \end{pmatrix}  &=
  \begin{pmatrix}
   w(t) \\ -z(t) + \rho \left(1 - z(t)^2\right)w(t)
 \end{pmatrix}
 \quad \forall t \in (0, T),
\end{align}
using the initial condition $z(0) = w(0) = 1$. Here, the design parameter $\rho>0$ determines the magnitude of the amplitude of the Van-der-Pol oscillator.

In order to set up an optimization problem, we choose a tracking type objective function that minimizes the discrepancy of the space-time averaged state $u(t,x):=\left(z(t), w(t), v(t,x)\right)^T$ to a prescribed value $a_{\text{target}}$:
\begin{align}\label{PinT:Numerics_objectivefunction}
    J = \frac 1 2\left(\frac 1 T \int_0^T \|u(t,\cdot)\|^2 \, \mathrm{d}t - a_{\text{target}} \right)^2 + \frac{\gamma}{2} \|\rho \|^2,
\end{align}
where a regularization term has been added with $\gamma = 10^{-6}$. The target value $a_{\text{target}}$ is computed in advance using a parameter $\rho_{\text{target}}=3$.

The time domain is discretized into $N=60000$ time steps on the finest time level with $t^i = i\Delta t$ for $i=1, \dots, N$ and $\Delta t = 0.0005$ giving the final time $t^N = 30$. The implicit Crank-Nicolson time-marching scheme approximates the transient term while the spatial derivatives  are approximated with a second order linear upwind scheme for the advection and central finite differences for the diffusive term. The spatial grid is chosen as $x_l = l\Delta x, l=1, \dots,L$ with $\Delta x = 0.01, L=100$.
% JBS Comment:  Can you say that this kind of long time domain and dt/dx ratio are "typical"?
% SG: Actually, I think ``typically'' a bigger time step dt would be chosen for this dx. So the ratio dt/dx is a bit small is this case. I had some problems with convergence of the nonlinear iterations on coarser time-grids when using bigger dt. Should I improve on the nonlinear solver and try to find a better combination?

At each time step, functional iterations are applied to solve the resulting nonlinear equations for $\vec u^i := (z^i, w^i, v^i_1, \ldots, v^1_L)\in \mathbb{R}^{2+L}$. These iterations are wrapped into the time-stepping function $\Phi^i$ required by the XBraid interface. The desired partial derivatives of $\Phi^i$ and $J$ that are needed in the adjoint interface are generated using the AD-Software CoDiPack \cite{CoDiPack,sagebaum2017high}.
The time-grid hierarchy for MGRIT uses a coarsening factor of $m=4$ and a maximum of three time-grid levels.
% Adding more levels generates coarse-grid time step sizes that are incompatible with the nonlinear time stepper.
% However, even with three levels, a reasonable speedup can be achieved which demonstrates the potential of the proposed method.

\begin{figure}
  \includegraphics{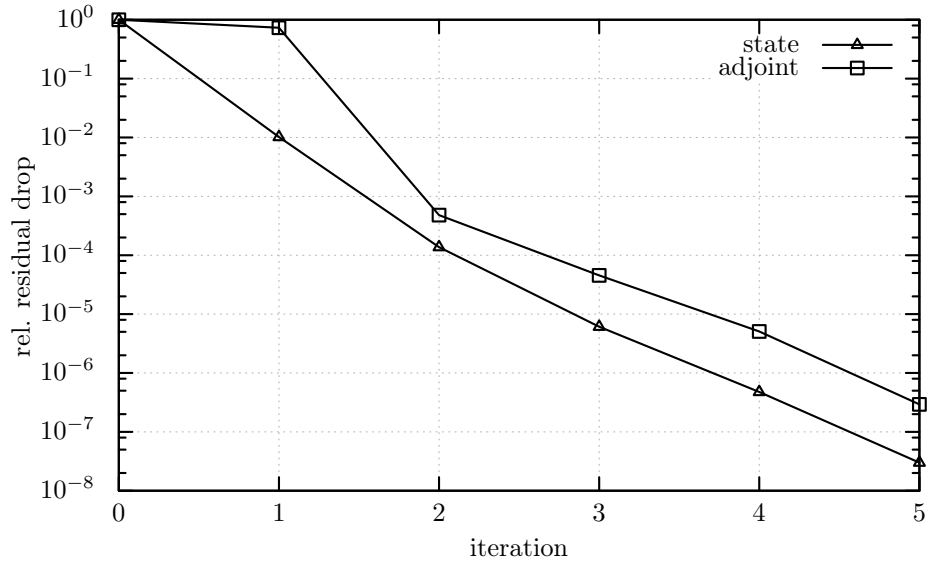}
  \caption{Relative drop of primal and adjoint residuals during piggyback iteration using the primal and adjoint XBraid solvers.}
  \label{fig:piggybackiter}
\end{figure}

\begin{figure}
  \includegraphics{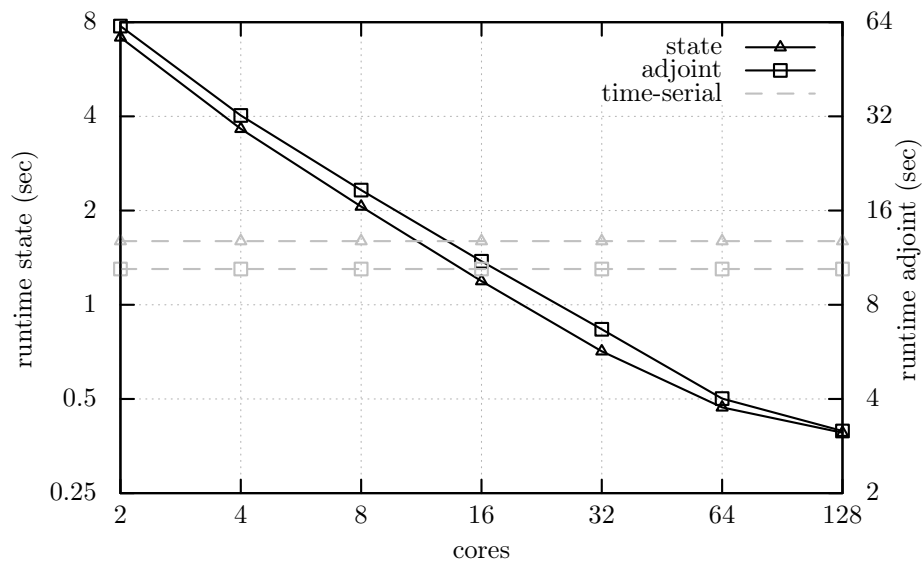}
  \caption{Strong scaling results for primal and adjoint time-parallelization during piggyback iterations. Grey dashed lines denote runtimes of the corresponding time-serial forward and backward state and adjoint time-stepping schemes.}
  \label{fig:piggyback_scaling}
\end{figure}

As a first step, the piggyback iteration has been implemented that solves the state and adjoint equations in time-parallel for a fixed design $\rho=2$. Figure \ref{fig:piggybackiter} shows the relative decrease in the state and adjoint residuals, $\|\vec u_{k+1} - \vec u_k\|_2$ and $\| \bar{\vec u}_{k+1} - \bar{\vec u}_k \|_2$. As expected, both residuals drop simultaneously while the adjoint iterates exhibit a certain time lag.
Strong scaling results are presented in Figure \ref{fig:piggyback_scaling}. It can be seen that speedup over the time-serial state and adjoint solvers is achieved if more that $16$ processors for time-parallelization are used. At $128$ processors in time, a speedup of about $3.15$ and $3.48$ for the state and the adjoint, respectively, has been achieved.

% \subsection{Parallel-in-time One-shot optimization}

As a next step, the parallel-in-time One-shot iteration integrates design updates into the piggyback iteration.
% It solves the inverse design problem for minimizing the tracking-type objective function where a target value has been computed in advance with $\rho_{\text{target}} = 3$.
The reduced gradient is preconditioned with a constant matrix $B_k = \theta \vec I$ using a constant step-size of $\theta = 0.9$ which was observed to yield sufficient descent on the Lagrange function.

\begin{figure}
  \center
  \includegraphics{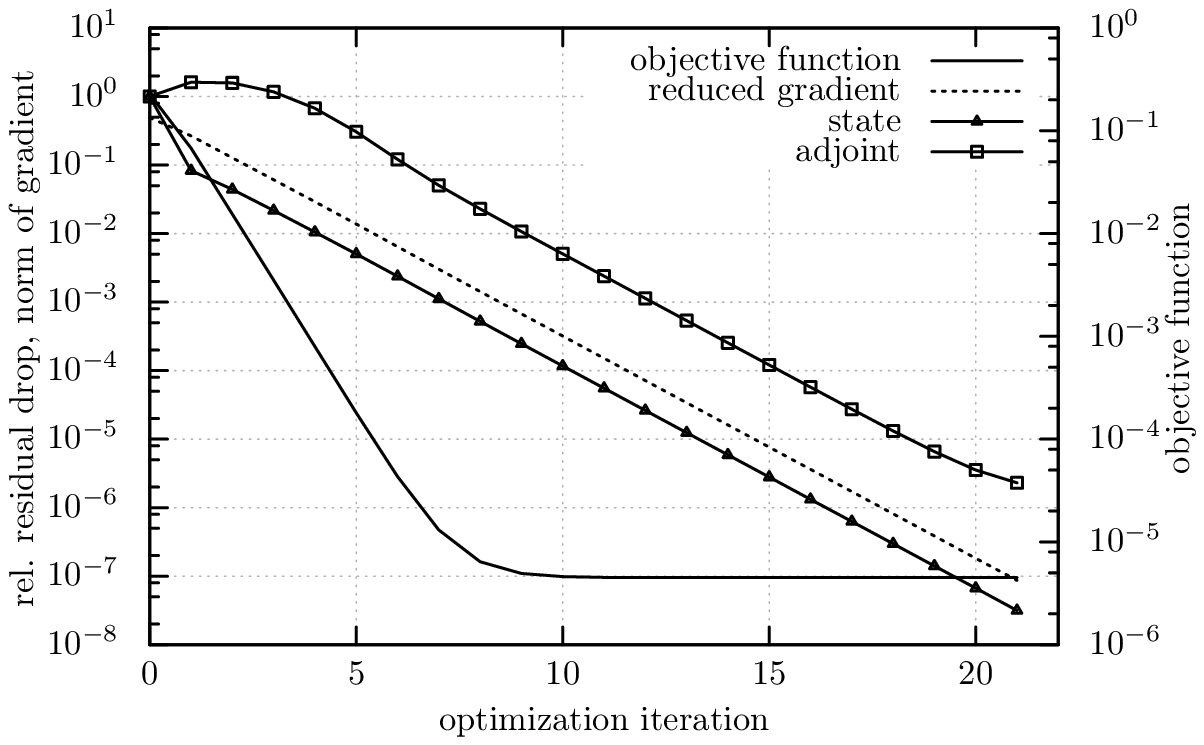}
  \caption{Optimization history of the parallel-in-time One-shot iteration.}
  \label{fig:PinT_oneshot_iteration}
\end{figure}

Figure \ref{fig:PinT_oneshot_iteration} plots the optimization history of the One-shot method. The reduced gradient drops simultaneously with the relative state and adjoint residuals. A total number of $22$ iterations is needed before the stopping criterion on the norm of the reduced gradient of $10^{-7}$ is reached. At that point, the state and adjoint residuals have been reduced sufficiently and the objective function has leveled out at the order of the regularization term.

In order to investigate the benefits of a parallel-in-time One-shot method, two reduced-space optimization methods are implemented. The first one is the conventional time-serial reduced-space optimization algorithm as in \eqref{reducedoptimiter1} - \eqref{reducedoptimiter3}. It employs the time-serial forward and backward time-marching scheme for the state and adjoint variables in each iteration of an outer optimization loop.
The second reduced-space optimization algorithm replaces the forward and backward time-marching schemes with the time-parallel state and adjoint MGRIT iterations. However, in contrast to the time-parallel One-shot method, it fully recovers a state and adjoint solution after each design update. This test case demonstrates the additional benefit of implementing the simultaneous optimization approach in a One-shot setting. Both reduced-space approaches utilize the same preconditioning matrix $B_k = \theta \vec I$, $\theta=0.9$ and the same stopping criterion of $10^{-7}$ on the norm of the reduced gradient.

\begin{figure}
  \label{oneshot_speedup}
  \includegraphics{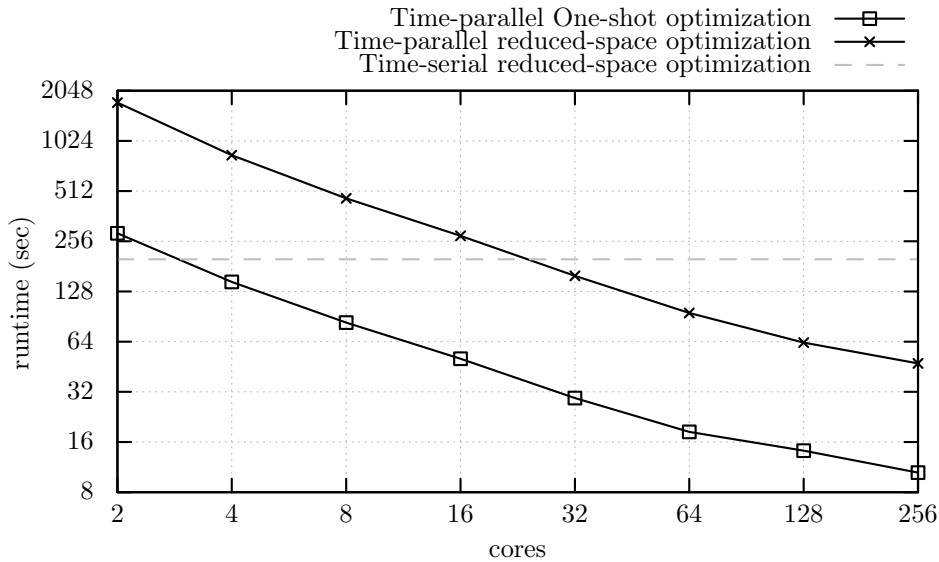}
  \caption{Strong scaling results for time-parallel One-shot optimization and time-parallel as well as time-serial reduced-space optimization.}
\end{figure}

Figure \ref{oneshot_speedup} compares runtimes of the three optimization approaches for increasing numbers of processors used for time-parallelization. The conventional time-serial optimization serves as a baseline as it is limited to one processor in time -- analogous to the situation where a spatially parallel code has reached its strong scaling limit. Both time-parallel approaches display the same slope for runtime reduction. However, the crossover point where the time-parallel approach yields a speedup over the time-serial approach is drastically reduced when using the simultaneous One-shot method. Here, speedup of the time-parallel One-shot method over the conventional time-serial optimization is achieved utilizing only $4$ processors for time-parallelization. Compared to the crossover point of the pure simulation in Figure \ref{fig:piggyback_scaling}, this further indicates that the time-parallel simultaneous optimization approach potentially provides speedup in cases where the pure time-parallel simulation approach fails to do so.

Runtime speedups when using $256$ processors for time parallelization are listed in Table \ref{tab:PinT:SpeedupOneshot}. The time-parallel reduced-space optimization achieves a speedup of 4.2 over the conventional method, which is at the same order as the reported speedup for the time-parallel state and adjoint computations. Beyond that, the One-shot framework yields an additional speedup factor of approximately $5$ over the time-parallel reduced-space method. Compared to the conventional time-serial optimization, a total speedup of $19$ is achieved for the parallel-in-time One-shot method.
% Even though more iterations are needed for the One-shot method, the runtime is reduced significantly because each iteration is computationally less expensive.

\begin{table}
  % \center
  \tbl{Speedup for the time-parallel One-shot method as well as the time-parallel reduced-space optimization compared to conventional time-serial reduced-space optimization.}
  {
  \begin{tabular}{ @ { } llrr @ { } }
% \hline\noalign{\smallskip}
    \toprule
            & cores & runtime & speedup \\
% \hline\noalign{\smallskip}
\colrule
    time-serial reduced-space optim.    &  1  & 199 sec  & 1.0 \\
    time-parallel reduced-space optim.  & 256 & 47 sec  & 4.2 \\
    time-parallel One-shot optim.       & 256 & 10 sec  & 19.0 \\
    \botrule
% \hline\noalign{\smallskip}
  \end{tabular}}
  \label{tab:PinT:SpeedupOneshot}
\end{table}

Table \ref{tab:overhead} compares the corresponding overhead factors of the optimization runtimes over the runtime of a pure time-serial simulation. The time-parallel One-shot algorithms reduces the total overhead significantly such that an optimization process can be achieved at the cost of only $6.3$ times the cost of the original time-serial simulation solver.

\begin{table}
  \center
  \tbl{Runtime overhead of optimization compared to pure time-serial simulation.}
  {
  \begin{tabular}{ @ { } llrr  @ { } }
    % \hline\noalign{\smallskip}
    \toprule
         & cores & runtime  &   Optim./Sim. \\
      % \hline\noalign{\smallskip}
      \colrule
      Simulation                         & 1   & 1.6 sec   &   1.0 \\
      Time-serial reduced-space optim.   & 1   & 199 sec   &  124.4 \\
      Time-parallel reduced-space optim. & 256 & 47 sec   &  29.4  \\
      Time-parallel One-shot optim.      & 256 & 10 sec   &  6.3  \\
      \botrule
    % \hline\noalign{\smallskip}
  \end{tabular}
  }\label{tab:overhead}
\end{table}

\section{Conclusion}
\label{sec:conclusion}
A non-intrusive framework for reducing the runtime of conventional gradient-based optimization algorithms for unsteady PDE-constrained optimization has been derived. The new framework applies parallel-in-time MGRIT iterations to the forward and backward time integration loops of the unsteady PDE and the adjoint equations, respectively, and integrates those into a simultaneous One-shot optimization approach. In this approach, design updates are employed after each state and adjoint update, based on a preconditioned inexact gradient.
Since the parallel-in-time approach offers the possibility to distribute workload onto multiple processors along the time domain, speedup over the conventional time-serial methods can be achieved through greater concurrency. Additionally, solving the unsteady dynamics iteratively with MGRIT enables the simultaneous One-shot optimization method to further reduce the runtime by applying design updates already at an early stage of the simulation process. The parallel-in-time One-shot method is non-intrusive to existing primal forward and adjoint backward time-stepping schemes such that transitioning from a conventional time-serial optimization algorithm requires only minimal additional coding. The benefit of the new framework has been demonstrated on an advection-dominated model problem. Here, a speedup of $19$ using $256$ processors in time has been achieved when compared to a conventional time-serial optimization algorithm.  This demonstrates the potential of the proposed method.

%
%%%%%%%%%%%%%%%%%%%%%%%%%%%%%%%%%%%%%%%%%%%%%%%%

\bibliographystyle{gOMS}
\bibliography{abbreviated}

\end{document}